\documentclass[12pt]{article}

\usepackage{amssymb,amsthm,amscd,array}

\let\ssection=\section
\renewcommand{\section}{\setcounter{equation}{0}\ssection}

\newcommand{\bbR}{\mathbb{R}}
\newcommand{\bbRP}{\mathbb{RP}}

\newcommand{\bbC}{\mathbb{C}}

\newcommand{\cD}{{\mathcal{D}}}

\newcommand{\cE}{{\mathcal{E}}}

\newcommand{\Diff}{\mathrm{Diff}}

\newcommand{\cF}{{\mathcal{F}}}

\newcommand{\rg}{\mathrm{g}}

\newcommand{\rR}{\mathrm{R}}

\newcommand{\SL}{\mathrm{SL}}
\newcommand{\Sl}{\mathrm{sl}}

\newcommand{\se}{\mathrm{e}}
\newcommand{\so}{\mathrm{o}}

\newcommand{\Vect}{\mathrm{Vect}}

\newcommand{\half}{\frac{1}{2}}

\def\mod#1{\left|{#1}\right|}

\begin{document}



\def\a{\alpha}
\def\b{\beta}
\def\d{\delta}
\def\g{\gamma}
\def\om{\omega}
\def\r{\rho}
\def\s{\sigma}
\def\vfi{\varphi}
\def\vr{\varrho}
\def\l{\lambda}
\def\m{\mu}
\def\n{\nu}
\def\implies{\Rightarrow}

\oddsidemargin .1truein
\newtheorem{thm}{Theorem}[section]
\newtheorem{lem}[thm]{Lemma}
\newtheorem{cor}[thm]{Corollary}
\newtheorem{pro}[thm]{Proposition}
\newtheorem{ex}[thm]{Example}
\newtheorem{rmk}[thm]{Remark}
\newtheorem{defi}[thm]{Definition}


\title{Conformally invariant differential operators on tensor densities}

\author{
V.~Ovsienko
\thanks{
CNRS, Centre de Physique Th\'eorique,
CPT-CNRS, Luminy Case 907,
F--13288 Marseille, Cedex~9, FRANCE;
mailto:ovsienko@cpt.univ-mrs.fr
}
\and
P. Redou
\thanks{Ecole Nationale d'Ingenieurs de Brest
Technopole Brest Iroise - Parvis Blaise Pascal
BP 30815 - 29608 Brest Cedex - FRANCE;
mailto:redou@enib.fr
}
}

\date{}

\maketitle

\bigskip

\thispagestyle{empty}

\begin{abstract}
Let $\cF_\l$ be the space of tensor densities on $\bbR^n$ of degree $\l$
(or, equivalently, of conformal densities of degree $-\l{}n$)
considered as a module over the Lie algebra $\so(p+1,q+1)$. We classify
$\so(p+1,q+1)$-invariant  bilinear differential operators from $\cF_\l\otimes\cF_\m$ 
to~$\cF_\n$. The classification of linear $\so(p+1,q+1)$-invariant
differential operators from
$\cF_\l$ to $\cF_\m$ already known in the literature (see \cite{ER,GJMS})
is obtained in a different manner.
\end{abstract}


\newpage

\section{Introduction}

\subsection{Linear and bilinear $\SL_2$-invariant operators}

In the one-dimensional case, the problem of classification of $\SL_2$-invariant
(bi)linear differential operators has already been treated in the classical literature.

Consider the action of $\SL(2,\bbR)$ on the space of functions in one variable, say on
$C^\infty(\bbRP^1)(\cong{}C^\infty(S^1))$ given by
\begin{equation}
f(x)\mapsto
f\left(\frac{ax+b}{cx+d}\right)(cx+d)^{-2\l}
\label{SLAction}
\end{equation}
depending on a parameter $\l\in\bbR$ (or $\bbC$). This $\SL(2,\bbR)$-module 
is called the space of tensor densities of degree $\l$ and denoted $\cF_\l$.
The classification of
$\SL(2,\bbR)$-invariant linear differential operators from $\cF_\l$ to $\cF_\m$
(i.e. of the operators commuting with the action (\ref{SLAction})) was obtained in classical
works on projective differential geometry. The result is as follows: there exists a unique
(up to a constant) $\SL(2,\bbR)$-invariant linear differential operator 
\begin{equation}
A_k: \cF_{\frac{1-k}{2}}
\to\cF_{\frac{1+k}{2}}
\label{InvOp1}
\end{equation}
where $k=0,1,2,\ldots$, and there are no $\SL(2,\bbR)$-invariant
operators from $\cF_\l$ to $\cF_\m$ for any other values of $\l$ and $\m$. Each operator
$A_k$ is of order $k$ and, for any choice of the coordinate $x$ such that the
$\SL(2,\bbR)$-action is as in (\ref{SLAction}),  is plainly given by the $k$-th order
derivative, $(A_kf)(x)=f^{(k)}(x)$.

The bilinear $\SL(2,\bbR)$-invariant differential operators from $\cF_\l\otimes\cF_\m$
to $\cF_\n$ were already classified by Gordan \cite{Gor}. For generic values of $\l$ and
$\m$ (more precisely, for $\l,\m\not=0,\half,1,\ldots$), there exists a unique (up to a
constant)
$\SL(2,\bbR)$-invariant bilinear differential operator 
\begin{equation}
B_k: \cF_\l\otimes\cF_\m
\to\cF_{\l+\m+k}
\label{BiInvOp1}
\end{equation}
where $k=0,1,2,\ldots$, and there are no $\SL(2,\bbR)$-invariant
operators for any other value of $\n$. This differential operator is given by the formula
\begin{equation}
B_k(f,g)=
\sum_{i+j=k}
(-1)^i{{2\m+k-1}\choose{i}}{{2\l+k-1}\choose{j}}f^{(i)}g^{(j)}
\label{Trans}
\end{equation}
and is called the transvectant.

\subsection{Multi-dimensional analogues}

In the multi-dimensional case, one has to distinguish the conformally flat case,
that can be reduced to $\bbR^n$ endowed with the standard $\so(p+1,q+1)$-action, where
$n=p+q$, and the curved (generic) case of an arbitrary manifold $M$ endowed with a
conformal structure.

In the conformally flat case, the analogues of the operators (\ref{InvOp1}) was classified
in \cite{ER}.  The result is as follows.

\begin{thm}(\cite{ER})
There exists a unique (up to a constant) $\so(p+1,q+1)$-invariant linear differential
operator 
\begin{equation}
A_{2k}: \cF_{\frac{n-2k}{2n}}
\to\cF_{\frac{n+2k}{2n}}
\label{InvOpn}
\end{equation}
where $k=0,1,2,\ldots$, and there are no $\so(p+1,q+1)$-invariant
operators from $\cF_\l$ to $\cF_\m$ for any other values of $\l$ and $\m\in\bbR$ (or
$\bbC$).
\label{LinThm}
\end{thm}
In the adopted coordinate system (corresponding
to the chosen conformally flat structure), the explicit expressions of the
operators $A_{2k}$  are
\begin{equation}
A_{2k}=
\Delta^k,
\qquad
\hbox{where}
\quad
\Delta = 
\rg^{ij}\frac{\partial}{\partial{}x^i}
\frac{\partial}{\partial{}x^j}
\label{InvOpn}
\end{equation}
In the generic (curved) case, the situation is much more complicated, see
\cite{GJMS, Gra}. (We also refer to \cite{Br} for a recent study of conformally invariant
differential operators on tensor fields, and a complete list of references.)

The purpose of this note is to extend the classical Gordan result to the multi-dimensional case.
We will classify $\so(p+1,q+1)$-invariant bilinear differential operators 
consider differential operators on tensor densities on $\bbR^n$, where $n=p+q$. 
The results can be also generalized for an arbitrary manifold $M$ endowed
with a conformally flat structure (e.g. upon a pseudo-Riemannian manifold $(M,\rg)$ with a
conformally flat metric). We will also give a simple direct proof of Theorem~\ref{InvOpn}.

We will not consider the curved case and state here a problem of existence of
bilinear conformally invariant differential operators for an arbitrary conformal structure.

\begin{rmk}
{\rm
There are other ways to generalize $\Sl_2$-symmetries in
the multi-dimen\-sional case. For instance, one can consider the $\Sl(n+1,\bbR)$-action
on~$\bbR^n$; this is related to projective differential geometry.
}
\end{rmk}

\section{Main result}

\subsection{Modules of differential operators}

Let $\cF_\l$ be the space of tensor densities of degree~$\l$ on $\bbR^n$, i.e. of smooth
section of the line bundle 
$\Delta_\l(\bbR^n)=
\mod{\Lambda^nT^*\bbR^n}^{\otimes\l}$ over $\bbR^n$.
We will be considering the space~$\cD_{\l,\m}$ of linear differential operators
from $\cF_\l$ to $\cF_\m$
and the space~$\cD_{\l,\m;\n}$ of bilinear differential operators
from $\cF_\l\otimes\cF_\m$ to $\cF_\n$.
These spaces of differential operators are naturally
$\Diff(\bbR^n)$- and $\Vect(\bbR^n)$-modules.
Note that the modules $\cD_{\l,\m}$ have been studied in a series of recent papers (see
\cite{DLO} and references therein).

Denote $\rg$ the standard quadratic form on $\bbR^n$ of signature $p-q$, where $p+q=n$.
The Lie algebra of infinitesimal conformal transformations is generated by the vector
fields
\begin{equation}
\matrix{
X_i  
&=&
\displaystyle \frac{\partial}{\partial x^i}\;\hfill\cr
\noalign{\smallskip}
X_{ij}
&=&
\displaystyle x_i\frac{\partial}{\partial x^j}-
x_j\frac{\partial}{\partial x^i}\;\hfill\cr
\noalign{\smallskip}
X_0  
&=&
\displaystyle x^i\frac{\partial}{\partial x^i}\;\hfill\cr
\noalign{\smallskip}
\bar X_i 
&=&
\displaystyle x_jx^j\frac{\partial}{\partial x^i}-
2x_ix^j\frac{\partial}{\partial x^j}\hfill\cr
}
\label{confGenerators}
\end{equation}
where $(x^1,\ldots,x^n)$ are coordinates on $\bbR^n$ and $x_i=\rg_{ij}x^j$, throughout this
paper, sum over repeated indices is understood.
Let us also consider the following Lie subalgebras
\begin{equation}
\so(p,q)\subset\se(p,q)\subset\so(p+1,q+1)\subset\Vect(\bbR^n)
\label{subalgebras}
\end{equation}
where $\so(p,q)$ is generated by the $X_{ij}$ and the Euclidean
sub\-algebra $\se(p,q)$ by $X_{ij}$ and $X_i$.

We will study the spaces $\cD_{\l,\m}$ and $\cD_{\l,\m;\n}$ as
$\so(p+1,q+1)$-modules and classify the differential operators commuting with the
$\so(p+1,q+1)$-action. 

It should be stressed that the classification of differential
operators invariant with respect to the Lie subalgebras $\so(p,q)$ is the
classical result of the Weyl theory of invariants \cite{Wey} (see also \cite{DLO} for the
case of $\se(p,q)$). We will use the Weyl classification in our work.

\begin{rmk}
{\rm
It is worth noticing that the conformal Lie algebra $\so(p+1,q+1)$ is maximal in the class
of finite-dimensional subalgebras of $\Vect(\bbR^n)$, that is, any bigger subalgebra of
$\Vect(\bbR^n)$ is infinite-dimensional (see \cite{BL} for a simple proof). 
}
\end{rmk}

\subsection{Introducing multi-dimensional transvectants}

The multi-dimensional analogues of the transvectants (\ref{Trans}) are described in the
following

\begin{thm}
For every $\l,\m\not=0,\frac{1}{n},\frac{2}{n},\ldots$, there exists a
unique (up to a constant) $\so(p+1,q+1)$-invariant bilinear differential operator 
\begin{equation}
B_{2k}: \cF_\l\otimes\cF_\m
\to\cF_{\l+\m+\frac{2k}{n}}
\label{BiInvOpn}
\end{equation}
where $k=0,1,2,\ldots$, and there are no $\so(p+1,q+1)$-invariant
operators from $\cF_\l\otimes\cF_\m$ to $\cF_\n$ for any other value of $\n$. 
\label{BiLinThm}
\end{thm}

The explicit formula for the operators $B_{2k}$ is complicated and is known only in some
particular cases.

\begin{rmk}
{\rm
The differential operators $A_{2k}$ and $B_{2k}$ are of order $2k$; comparing with the
one-dimensional case, one has twice less invariant differential operators. Note that if one
takes semi-integer $k$ in  (\ref{InvOpn}), then the corresponding operator is
pseudo-differential.
}
\end{rmk}

It would be interesting to obtain a complete classification of
$\so(p+1,q+1)$-invariant bilinear differential operators (see \cite{Gar} for the
one-dimensional case).

\section{Proof of Theorems \ref{LinThm} and \ref{BiLinThm}}

We will start the proof with classical results of the theory of invariants and describe the
differential operators invariant with respect to the action of the Lie algebra~$\se(p,q)$.
We refer \cite{Wey} as a classical source and \cite{DLO} for the description of the
Euclidean invariants.

\subsection{Euclidean invariants}\label{EucInv}

Using the standard affine connection on $\bbR^n$, one identifies the space of linear
differential operators on $\bbR^n$ with the corresponding space of symbols, i.e., with the
space of smooth functions on
$T^*\bbR^n\cong\bbR^n\oplus(\bbR^n)^*$ polynomial on~$(\bbR^n)^*$.
This identification is an isomorphism of modules over the algebra of affine
transformations and allows us to apply the theory of invariants.

Moreover, choosing a
(dense) subspace of symbols which are also polynomials on the first summand, one
reduces the classification of
$\se(p,q)$-invariant differential operators from
$\cD_{\l,\m}$ to the classification of $\se(p,q)$-invariant polynomials in the space
$\bbC[x^1,\ldots,x^n,\xi_1,\ldots,\xi_n]$, where $(\xi_1,\ldots,\xi_n)$ are the coordinates
on$(\bbR^n)^*$ dual to $(x^1,\ldots,x^n)$.

Consider first invariants with respect to $\so(p,q)\subset\se(p,q)$. It is well-known
(see~\cite{Wey}) that the algebra of $\so(p,q)$-invariant polynomials is generated by three
elements
\begin{equation}
\rR_{xx}=\rg_{ij}\,x^ix^j,
\quad
\rR_{x\xi}=x^i\xi_i,
\quad
\rR_{\xi\xi}=\rg^{ij}\,\xi_i\xi_j
\label{InvPol}
\end{equation}
Second, taking into account the invariance with respect to translations in $\se(p,q)$, any
$\se(p,q)$-invariant polynomial $P(x,\xi)$ satisfies $\partial{P}/\partial{}x^i=0$.
The only remaining generator is $\rR_{\xi\xi}$ and, therefore, $\se(p,q)$-invariant
linear differential operators from $\cF_\l$ to $\cF_\m$ are linear combinations of
operators (\ref{InvOpn}).

Note that the obtained result is, of course, independent on $\l$ and $\m$ since the degree
of tensor densities does not intervene in the $\se(p,q)$-action.

\subsection{Proof of Theorem \ref{LinThm}}\label{SecProofLin}

We must check now for which values of $\l$ and $\m$ the operators (\ref{InvOpn}) from
$\cF_\l$ to $\cF_\m$ are invariant with respect to the action of the full conformal
algebra $\so(p+1,q+1)$. 

By definition, the action of a vector field $X$ on an element
$A\in\cD_{\l,\m}$ is given by
\begin{equation}
L_X^{\l,\m}(A)=
L_X^\m\circ{}A-A\circ{}L_X^\l
\label{DefAction}
\end{equation}
where $L_X^\l$ is the operator of Lie derivative of $\l$-density, namely
\begin{equation}
L_X^\l=X^i\frac{\partial}{\partial{}x^i}+\l\,\partial_iX^i
\label{LieDer}
\end{equation}
for any coordinate system.

Consider the action of the generator $X_0$ in (\ref{confGenerators}) on the operator
$A=\sum_{k\geq0}c_k\rR_{\xi\xi}^k$. Using the preceding expressions, one readily gets
\begin{equation}
L_{X_0}^{\l,\m}(A)=
\sum_{i\geq0}\left(
n\d-2k
\right)
c_k\rR_{\xi\xi}^k
\label{DefAction}
\end{equation}
where $\d=\m-\l$.
Thus, the invariance condition $L_{X_0}^{\l,\m}(A)=0$ is satisfied if and only if for each $k$
in the above sum either $c_k=0$ or $\d=\frac{2k}{n}$; and one gets the values of the
shift $\d$ in accordance with~(\ref{InvOpn}).

\goodbreak

Consider, at last, the action of the generators $\bar{X}_i$ (with $i=1,\ldots,n$). After the
identification of the differential operators with polynomials one has the following
result from  \cite{DLO}.

\begin{pro}
The action of the generator $\bar{X}_i$ on $\cD_{\l,\m}$ is as follows
\begin{equation}
L_{\bar X_i}^{\l,\m}=
L_{\bar X_i}^\d-
\xi_iT + 2(\cE+n\l)\,\partial_{\xi^i}
\label{InversionAction}
\end{equation}
where 
\begin{equation}
L^\d_{\bar{X}_i}
=
x_jx^j\partial_i - 2x_ix^j\partial_j
-2(\xi_i x_j - \xi_j x_i)\partial_{\xi_j}
+2\xi_j x^j\partial_{\xi^i}
-2n\d{}x_i
\label{InversionLift}
\end{equation}
is the cotangent lift, and where $T=\partial_{\xi^j}\partial_{\xi_j}$ is the trace and 
$\cE=\xi_j\partial_{\xi_j}$ the Euler operator.
\label{FirstPro}
\end{pro}

Applying $L_{\bar X_i}^{\l,\m}$ to the operator $\rR^k_{\xi\xi}$ one then obtains
$$
\begin{array}{rcl}
L^{\l,\m}_{\bar{X}_i}(\rR_{\xi\xi}^k)
&=&2(2k-n\d)x_i\rR_{\xi\xi}^k\\[6pt]
&&+2k(n(2\l-1)+2k)\xi_i\rR_{\xi\xi}^{k-1}
\end{array}
$$
The first term in this expression vanishes for $2k-n\d$=0, this condition is precisely the
preceding one; the second term vanishes if and only if
$$
\l=\frac{n-2k}{2n}.
$$
Hence the result.

\subsection{Proof of Theorem \ref{BiLinThm}}

As in Section \ref{EucInv}, let us first consider the operators invariant with respect to
the  Lie algebra $\se(p,q)$. Again, identifying the bilinear differential operators with their
symbols, one is led to study the algebra of $\se(p,q)$-invariant polynomials in the space
$\bbC[x^1,\ldots,x^n,\xi_1,\ldots,\xi_n,\eta_1,\ldots,\eta_n]$. The Weyl invariant theory
just applied guarantees that there are three generators:
\begin{equation}
\rR_{\xi\xi}=\rg^{ij}\,\xi_i\xi_j,
\quad
\rR_{\xi\eta}=\rg^{ij}\,\xi_i\eta_j,
\quad
\rR_{\eta\eta}=\rg^{ij}\,\eta_i\eta_j
\label{BiEucInv}
\end{equation}
Any $\se(p,q)$-invariant bilinear differential operator is then of the form
$$
B=\sum_{r,s,t\geq0}
c_{rst}\rR^{r,s,t}
$$
where, to simplify the notations, we put
$\rR^{r,s,t}=\rR_{\xi\xi}^r
\rR_{\xi\eta}^s
\rR_{\eta\eta}^t.
$

\goodbreak

The action of a vector field $X$ on a bilinear operator 
$B:\cF_\l\otimes\cF_\m\to\cF_\n$ is defined as follows
\begin{equation}
(L_{X}^{\l,\m;\n}B)(f,g)=
L_X^\n{}B(f,g)-B(L_X^\l{}f,g)-B(f,L_X^\m{}g)
\label{BiAction}
\end{equation}

Let us apply the generator $X_0$ to the operator $B$, one has
\begin{equation}
L_{X_0}^{\l,\m;\n}(B)=
\sum_{r,s,t\geq0}\left(
n\d-2(r+s+t)
\right)
c_{rst}\rR^{r,s,t}
\label{ZeroAction}
\end{equation}
where $\d=\n-\m-\l$. The equation $L_{X_0}^{\l,\m}(B)=0$ leads to the homogeneity
condition 
\begin{equation}
\d=\frac{2(r+s+t)}{n}
\label{Homogen}
\end{equation}
The general expression for $B$ retains the form
\begin{equation}
B_{2k}=\sum_{r+s+t=k}
c_{rst}\rR^{r,s,t}
\label{BiInvPol}
\end{equation}
The operator $B$ is of order $2k$ and $\n=\l+\m+2k/n$.

Now, to determine the coefficients $c_{rst}$, one has to apply the generators $\bar{X}_i$.
One has the following analog of Proposition \ref{FirstPro}.
\begin{pro}
The action of the generator $\bar{X}_i$ on $\cD_{\l,\m;\n}$ is given by
\begin{equation}
L_{\bar X_i}^{\l,\m;\n}=
L_{\bar X_i}^\d-
\xi_iT_{\xi}-\eta_iT_\eta + 2\left(
(\cE_\xi+n\l)\,\partial_{\xi^i}
+(\cE_\eta+n\m)\,\partial_{\eta^i}
\right)
\label{InversionAction}
\end{equation}
where 
\begin{equation}
\begin{array}{rcl}
L^\d_{\bar{X}_i}
&=&
x_jx^j\partial_i - 2x_ix^j\partial_j-2n\d{}x_i\\[6pt]
&&-2\left(
(\xi_i x_j - \xi_j x_i)\partial_{\xi_j}
+(\eta_i x_j - \eta_j x_i)\partial_{\eta_j}
\right)\\[6pt]
&&+2\left(
\xi_j x^j\partial_{\xi^i}
+\eta_j x^j\partial_{\eta^i}
\right)
\end{array}
\label{InversionDoubleLift}
\end{equation}
is just the natural lift of $\bar{X}_i$ to $T^*\bbR^n\oplus{}T^*\bbR^n$.
\end{pro}
Applying $L_{\bar X_i}^{\l,\m;\n}$ to each monomial,  
$\rR^{r,s,t}=\rR_{\xi\xi}^r
\rR_{\xi\eta}^s
\rR_{\eta\eta}^t$, in the operator $B_{2k}$, one immediately gets
\begin{equation}
\begin{array}{rcl}
L_{\bar X_i}^{\l,\m;\n}(\rR^{r,s,t})
&=&
2(2k-n\d)x_i\rR^{r,s,t}\\[6pt]
&&
\begin{array}{l}
+\Big(
2r(2r+n(2\l-1))\rR^{r-1,s,t}
-s(s-1)\rR^{r,s-2,t+1}
\\[4pt]
\;+2s(s+2t+n\m-1)\rR^{r,s-1,t}
\Big)\xi_i
\end{array}
\\[6pt]
&&
\begin{array}{l}
+\Big(
2t(2t+n(2\m-1))\rR^{r,s,t-1}
-s(s-1)\rR^{r+1,s-2,t}
\\[4pt]
\;+2s(s+2r+n\l-1)\rR^{r,s-1,t}
\Big)\eta_i
\end{array}
\end{array}
\label{InversionDouble Act}
\end{equation}
At last, applying $L_{\bar X_i}^{\l,\m;\n}$ to the operator $B_{2k}$ written in the form
(\ref{BiInvPol}) and collecting the terms, one readily gets the following recurrent system
of two linear equations
\begin{equation}
\begin{array}{rcl}
2(r+1)\left(
2(r+1)+n(2\l-1)
\right)c_{r+1,s,t}&&\\[4pt]
-(s+2)(s+1)\,{}c_{r,s+2,t-1}&&\\[4pt]
+2(s+1)(s+2t+n\m)\,{}c_{r,s+1,t}&=&0
\end{array}
\label{System1}
\end{equation}
and
\begin{equation}
\begin{array}{rcl}
2(t+1)\left(
2(t+1)+n(2\m-1)
\right)c_{r,s,t+1}&&\\[4pt]
-(s+2)(s+1)\,{}c_{r-1,s+2,t}&&\\[4pt]
+2(s+1)(s+2r+n\l)\,{}c_{r,s+1,t}&=&0
\end{array}
\label{System2}
\end{equation}
for the coefficients. For $\l,\m\not=0,\frac{1}{n},\frac{2}{n},\ldots$ this system
has a unique (up to a constant) solution. Indeed, choosing $c_{0,k,0}$ as a parameter, one
uses the first equation to express $c_{r+1,s,t}$ from $c_{r,*,*}$ and the second
one to express $c_{r,s,t+1}$ from $c_{*,*,t}$.

Theorem \ref{BiLinThm} is proved.

\medskip

The differential operators (\ref{InvOp1}) and (\ref{BiInvOp1}) play important r\^ole in the
theory of modular functions (see \cite{ran,coh}), in projective differential geometry
(see \cite{Wil, Pee}) and in the representation theory of $\SL(2,\bbR)$. The transvectants
have been recently used in
\cite{cmz,Ovs} to construct $\SL(2,\bbR)$-invariant star-products on $T^*S^1$. We plan to
discuss the relation of the conformally invariant operators described in this note to the
representation theory in a subsequent paper.

\section{Examples}

Let us give here explicit formul{\ae} for the bilinear differential operators 
(\ref{BiInvOpn}) with $k=1,2$.
Using (\ref{System1}) and (\ref{System2}) one has:
$$
\begin{array}{rcl}
B_{2}
&=&
n\m(2+n(2\m-1))\rR_{\xi\xi}-(2+n(2\m-1))(2+n(2\l-1))\rR_{\xi\eta}\\[6pt]
&&
\begin{array}{l}
+n\l(2+n(2\l-1))\rR_{\eta\eta}
\label{secondorder}
\end{array}
\end{array}
$$
and 
$$
\begin{array}{rcl}
B_{4}
&=&
-(2+n(2\l-1))(2+n(2\m-1))(4+n(2\l-1))(4+n(2\m-1))\rR_{\xi\eta}^2\\[6pt]
&&
\begin{array}{l}
+2(1+n\m)(2+n(2\m-1))(4+n(2\l-1))(4+n(2\m-1))\rR_{\xi\xi}\rR_{\xi\eta}
\\[4pt]
+2(1+n\l)(2+n(2\l-1))(4+n(2\l-1))(4+n(2\m-1))\rR_{\xi\eta}\rR_{\eta\eta}
\end{array}
\\[6pt]
&&
\begin{array}{l}
-\frac{1}{2}\Big((2+n(2\l-1))+2(1+n\m)(2+n\l)+(2+n(2\m-1))
\\[4pt]
+2(1+n\l)(2+n\m)\Big)
(4+n(2\l-1))(4+n(2\m-1))\rR_{\xi\xi}\rR_{\eta\eta}
\\[4pt]
-(1+n\m)(2+n(2\m-1))n\m(4+n(2\m-1))\rR_{\xi\xi}^2
\\[4pt]
-(1+n\l)(2+n(2\l-1))n\l(4+n(2\l-1))\rR_{\eta\eta}^2.
\end{array}
\end{array}
$$
The expressions of further orders operators are much more 
complicated, and we do not have an explicit general formula, except for some 
particular coefficients $c_{r,s,t}$.
A direct computation using equations
 leads to  
$$
c_{i,k-i,0}=
(-1)^i{{k}\choose{i}}\frac{(k-1+n\m)(k-2+n\m)\cdots(k-i+n\m)}{(2+n(2\l-1))
(4+n(2\l-1))\cdots(2i+n(2\l-1))}c_{0,k,0}
$$
and
$$
c_{0,k-i,i}=
(-1)^i{{k}\choose{i}}\frac{(k-1+n\l)(k-2+n\l)\cdots(k-i+n\l)}{(2+n(2\m-1))
(4+n(2\m-1))\cdots(2i+n(2\m-1))}c_{0,k,0}.
$$
\medskip

\textbf{Acknowledgments:} We are indebted to C.~Duval, A.~A.~Kirillov,
Y.~Kosmann-Schwarzbach and P.~Lecomte for enlightening discussions.



\begin{thebibliography}{99}

\bibitem{BL}
F.~Boniver and P.~B.~A.~Lecomte,
{\it A remark about the Lie algebra of infinitesimal conformal transformations of the
Euclidian space}, Bull. London Math. Soc. {\bf  32:3} (2000) 263--266.

\bibitem{Br}
T. Branson, {\it Second order conformal covariants},
Proc. Amer. Math. Soc. {\bf  126:4} (1998) 1031--1042.

\bibitem{coh}
H. Cohen, {\it Sums involving the values at negative
integers of $L$ functions of quadratic characters},
Math. Ann. {\bf  217} (1975) 181--194.

\bibitem{cmz} P. Cohen, Yu. Manin \& D. Zagier,
{\it Automorphic pseudodifferential operators}, Algebraic
aspects of integrable systems, 17--47, Progr. Nonlinear Differential
Equations Appl., 26, Birkh\"auser Boston, Boston, MA, 1997.

\bibitem{DLO} C. Duval, P. Lecomte \& V. Ovsienko,
{\em 	Conformally equivariant quantization: existence and uniqueness},
Ann. Inst. Fourier. {\bf  49:6} (1999) 1999--2029.

\bibitem{ER}
M.G. Eastwood \& J.W. Rice, {\it Conformally invariant differential
operators on Minkowski space and their curved analogues},
Comm. Math. Phys. {\bf 109:2} (1987), 207--228.

\bibitem{Gar}
D. Garajeu, {\it Conformally and projective covariant differential operators},
Lett. Math. Phys. {\bf 47:4} (1999) 293--306.

\bibitem{Gor} P. Gordan, {\em Invariantentheorie}, Teubner, Leipzig, 1887.

\bibitem{GJMS} 
C.R. Graham, R. Jenne, L.J. Mason \&G.A. Sparling,
{\it Conformally invariant powers of the Laplacian. I. Existence},
J. London Math. Soc. (2) {\bf 46:3}
(1992), 557--565.

\bibitem{Gra} 
C.R. Graham,
{\it Conformally invariant powers of the Laplacian. II. Nonexistence},
J. London Math. Soc. (2)  {\bf 46:3} (1992),  566--576. 

\bibitem{Pee} 
S. Janson \& J. Peetre, {\it A new generalization of
Hankel operators (the case of higher weights)}, Math. Nachr. 132 (1987)
313--328.

\bibitem{Kos} 
Y. Kosmann,
{\it Sur les degr\'es conformes des op\'erateurs diff\'erentiels},
C. R. Acad. Sc. Paris, {\bf 280} (1975), 229--232.

\bibitem{Ovs}
 V. Ovsienko,
{\em Exotic Deformation Quantization},
 J. Diff. Geom.,
{\bf 45:2} (1997) 390--406.

\bibitem{PR}
R.~Penrose and W.~Rindler,
{\sl Spinors and space-time,
Vol.~2,
Spinor and twistor methods in space-time geometry},
Cambridge University Press, 
1986.

\bibitem{ran} R. A. Rankin, {\it The construction of
automorphic forms from the derivatives of a given
form}, J. Indian Math. Soc. 20 (1956) 103-116.

\bibitem{Wey}
H.~Weyl,
{\sl The Classical Groups}, 
Princeton University Press, 1946.

\bibitem{Wil} 
E.~J.~Wilczynski, 
{\sl Projective differential geometry
of curves and ruled surfaces}, 
Leipzig, Teubner, 1906.

\end{thebibliography}
\end{document}